\documentclass[11pt]{article}
\usepackage{amssymb,amsfonts,amsmath,latexsym,epsf,tikz,url}
\usepackage[usenames,dvipsnames]{pstricks}
\usepackage{pstricks-add}
\usepackage{epsfig}
\usepackage{pst-grad} 
\usepackage{pst-plot} 
\usepackage[space]{grffile} 
\usepackage{etoolbox} 
\makeatletter 
\patchcmd\Gread@eps{\@inputcheck#1 }{\@inputcheck"#1"\relax}{}{}
\makeatother

\newtheorem{theorem}{Theorem}[section]
\newtheorem{proposition}[theorem]{Proposition}

\newtheorem{corollary}[theorem]{Corollary}
\newtheorem{lemma}[theorem]{Lemma}

\newtheorem{problem}[theorem]{Problem}
\newcommand{\proof}{\noindent{\bf Proof.\ }}
\newcommand{\qed}{\hfill $\square$\medskip}
\newcommand{\gsp}{\gamma_{\rm sp}}

\begin{document}

\title{Computational complexity aspects of super domination}

\author{
Csilla Bujt\'as$^{1,4}$ 
\and
Nima Ghanbari$^{2}$ 
\and
Sandi Klav\v{z}ar$^{1,3,4}$
}

\date{\today}

\maketitle

\begin{center}
$^1$ Faculty of Mathematics and Physics, University of Ljubljana, Slovenia\\
{\tt csilla.bujtas@fmf.uni-lj.si}\\
{\tt sandi.klavzar@fmf.uni-lj.si}
\medskip

$^{2}$ Department of Informatics, University of Bergen, \\
P.O.\ Box 7803, 5020 Bergen, Norway\\
{\tt Nima.Ghanbari@uib.no}
\medskip

$^3$ Faculty of Natural Sciences and Mathematics, University of Maribor, Slovenia\\
\medskip

$^4$ Institute of Mathematics, Physics and Mechanics, Ljubljana, Slovenia\\
\medskip
\end{center}


\begin{abstract}
Let $G$ be a graph. A dominating set $D\subseteq V(G)$ is a super dominating set if for every vertex  $x\in V(G) \setminus D$ there exists $y\in D$ such that $N_G(y)\cap (V(G)\setminus D)) = \{x\}$. The cardinality of a smallest super dominating set of $G$ is the super domination number of $G$. An exact formula for the super domination number of a tree $T$ is obtained and demonstrated that a smallest super dominating set of $T$ can be computed in linear time. It is proved that it is NP-complete to decide whether the super domination number of a graph $G$ is at most a given integer if $G$ is a bipartite graph of girth at least $8$. The super domination number is determined for all $k$-subdivisions of graphs. Interestingly, in half of the cases the exact value can be efficiently computed from the obtained formulas, while in the other cases the computation is hard. While obtaining these formulas, II-matching numbers are introduced and proved that they are computationally hard to determine. 
\end{abstract}


\medskip
\noindent{\bf AMS Subj.\ Class.:} 05C69, 68Q25

\section{Introduction}
 
 Let $G = (V(G), E(G))$ be a graph. Then $D\subseteq V(G)$ is a  \textit{dominating set} if every vertex in $\overline{D}= V(G) \setminus D$ is adjacent to at least one vertex in $D$. The  \textit{domination number} $\gamma(G)$ of $G$ is the minimum cardinality of a dominating set of $G$. Graph domination theory has been extensively researched so far. To capture the current state of the field, we refer the reader to two recent edited books~\cite{haynes-2020, haynes-2021}. 

Many variations of the domination have been introduced, some of which are significant and important (such as total domination and connected domination), while others are of only sporadic importance. In our view, the group of significant domination concepts includes super domination which was introduced in 2015 by Lema\'nska, Swaminathan, Venkatakrishnan, and Zuazua~\cite{Lemans}.  

It is a classically known fact, that a dominating set $D$ of $G$ is minimal if and only if every vertex from $D$ has a private neighbor in $\overline{D} = V(G)\setminus D$. As a dual concept one says that a dominating set $D$ of $G$ is a \textit{super dominating set} of $G$, if for every vertex  $x\in \overline{D}$ there exists $y\in D$ such that $N_G(y) \cap \overline{D} = \{x\}$. (As usual, $N_G(v)$ stands for the open neighborhood of $v$ in $G$ and $N_G[v]$ for the closed neighborhood of $v$ in $G$.) 
The cardinality of a smallest super dominating set of $G$ is the \textit{super domination number} $\gsp(G)$ of $G$. A super dominating set of cardinality $\gsp(G)$ is briefly called a {\em $\gsp$-set}. The initial study of the concept has been followed by several sequels, of which the reader is referred to~\cite{Dett, Kle, Kri, sen-2022, Zhu}. 

The {\em subdivision} of a graph $G$ is the graph obtained from $G$ by replacing each edge with a disjoint path of length $2$ and is denoted by $S(G)$. More generally, if $k\ge 1$, then the graph $S_k(G)$ is obtained from $G$ by replacing each edge with a disjoint path of length $k+1$, that is, subdividing each of its edges $k$ times. Clearly, $S(G)=S_1(G)$. Some authors use the term {\em complete $k$-subdivision} for what we call $k$-subdivision, but since all our subdivisions are complete, we  simplify the terminology. The concept of $k$-subdivisions is ubiquitous in graph theory, here we cite its presence in graph colorings~\cite{KSUB, wang-2018}, spectral graph theory~\cite{barrett-2014}, structural graph theory~\cite{lozin-2022, zhou-2017}, and chemical graph theory~\cite{arock-2021, KlKnMa-2021}. 

In this article we discuss various aspects of super domination, which are in one way or another intertwined with the computational complexity of the problem of determining the super domination number. In the next section we recall some definitions and known results, and state a useful characterization of super dominating sets. Super domination has already been considered on trees from several perspectives, see~\cite{Kri, Lemans, sen-2022, Zhu}. By now, only sharp upper and lower bounds have been obtained. In Section~\ref{sec:trees} we fill this gap by providing an exact formula for the super domination number of a tree. Moreover, we demonstrate that if $T$ is a tree, then $\gsp(T)$ as well as a $\gsp$-set of $T$ can be computed in linear time. On the negative side, in Section~\ref{sec:bipartite}, we prove that it is NP-complete to decide whether $ \gsp(G) \leq k$  holds if $G$ is a bipartite graph of girth $g(G) \ge 8$ and the positive integer $k$  is part of the input. In our longest part of the paper, Section~\ref{sec:subdivisions}, we consider $k$-subdivisions of arbitrary graphs. Depending on $k \bmod 4$, four closed formulas for $\gsp(S_k(G))$ are proved. When $k \bmod 4 \in \{1,3\}$, the corresponding formulas depend only of $k$, the size of $G$, and a simple condition on the cycles of $G$. Note that in these two cases $S_k(G)$ is bipartite. On the other hand, if $k \bmod 4 = 0$, then $\gsp(S_k(G))$ is a function of $\gsp(G)$ also, and if $k \bmod 4 = 2$, then $\gsp(S_k(G))$ depends also on the cardinality of a largest matching that admits a partition into two induced matchings. We name such matchings as II-matchings and prove that it is NP-hard to compute the maximum size of such matchings. It follows that for each even $k$, it is also NP-hard to determine $\gsp(S_k(G))$.

\section{Preliminaries}
\label{sec:prelim}

The order and the size of a graph $G$ will be denoted by $n(G)$ and $m(G)$, respectively. If $D$ is a super dominating set of $G$ and if for a vertex $x\in \overline{D}$ the vertex $y\in D$ has the property $N_G(y) \cap \overline{D} = \{x\}$, then we will say that $x$ is {\em super dominated} by $y$.  

Let $G$ be a graph. Then the independence number of $G$ will be denoted by $\alpha(G)$, the matching number of $G$ by $\alpha'(G)$, and the vertex cover number of $G$ by $\beta(G)$. A set $X\subseteq V(G)$ is a {\em $2$-packing} of $G$ if $d_G(u,v)\ge 3$ holds for each pair of vertices $u, v\in X$. In other words, each pair of vertices of $X$ has disjoint closed neighborhoods. The cardinality of a smallest $2$-packing of $G$ will be denoted by $\rho(G)$. 

The path in $S_k(G)$ obtained by $k$ times subdividing an edge $uv\in E(G)$ will be denoted by $P_{uv}$ and addressed to as a {\em super edge}. The vertices of $P_{uv}$ will be denoted by $u,(uv)_1, \ldots, (uv)_{k},v$.  Note that $n(S_k(G)) = n(G) + k\cdot m(G)$ and $m(S_k(G)) = (k+1)m(G)$. We say that a graph is a {\em $k$-subdivision graph} if it can be obtained as a $k$-subdivision of some graph. 

We next recall a few results on the super domination number needed later on. 

\begin{theorem} {\rm \cite{Lemans}}
\label{thm:general-bounds}
If $G$ is a graph without isolated vertices, then,
$$1\leq \gamma(G) \leq \frac{n}{2} \leq \gsp(G) \leq n(G)-1\,.$$
\end{theorem}

\begin{theorem} {\rm \cite[Theorem 3, Corollary 2]{Kle}} 
\label{thm:alpha'-rho}
If $G$ is a graph with $n(G)\ge 2$, then 
$$n(G) - \alpha'(G) \le \gsp(G) \le n(G) - \rho(G).$$ 
\end{theorem}

Since $n(G) = \alpha(G) + \beta(G)$ in general and $\alpha'(G) = \beta(G)$ when $G$ is bipartite, Theorem~\ref{thm:alpha'-rho} implies that if $G$ is bipartite, then  $\gsp(G) \ge \alpha(G)$. 

\begin{theorem} {\rm \cite{Lemans}}
\label{thm:paths-cycles-stars}
The following exact values are valid.
\begin{itemize}
\item[(i)]
If $n\ge 2$, then $\gsp(P_n)=\lceil \frac{n}{2} \rceil$.
\item[(ii)]
If $n\ge 3$, then 
\begin{displaymath}
\gsp(C_n)= \left\{ \begin{array}{ll}
\lceil\frac{n+1}{2}\rceil; & n \equiv 2 \bmod 4, \\
\\
\lceil\frac{n}{2}\rceil; & \mbox{otherwise}.
\end{array} \right.
\end{displaymath}
\item[(iii)]
If $n\ge 2$, then $\gsp(K_{1,n})=n$.
\end{itemize}
	\end{theorem}

Let $G$ be a graph, and let  $D$ be a super dominating set of $G$. For each $u\in \overline{D}$, select an arbitrary vertex $u'\in D$ such that $u$ is the unique neighbor of $u'$ in $\overline{D}$, that is, $N(u') \cap \overline{D} = \{u\}$. Then we say that the set 
$$D^\ast = \{ u'\in D:\ u\in \overline {D}\}$$
is a {\em core of} $D$. By this definition, there exists a matching between $\overline{D}$ and $D^\ast$ that covers $\overline{D} \cup D^\ast$. Moreover, the following result holds which seems of independent interest. Before stating it, we introduce a notation. For two disjoint vertex sets $A, B \subseteq V(G)$ let $E_G[A,B]$ denote the set of all edges between $A$ and $B$ in $G$.

\begin{lemma} \label{lem:matching}
Let $A$ and $B$ be two disjoint vertex sets of a graph $G$. Then $D = \overline{A}$ is a super dominating set and $B$ is a core of $D$ if and only if $E_G[A,B]$ is a matching that covers all vertices in $A \cup B$.
\end{lemma}

\proof 
Suppose that $D = \overline{A}$ is a super dominating set and let $B$ be its core. 
By definition of the core, every vertex from $B$ has exactly one neighbor from $A$, so that $|N_G(x) \cap B|=1$ if $x\in A$ and $|N_G(y)\cap A|=1$ whenever $y \in B$. It shows that $E_G[A,B]$ is a matching that covers $A \cup B$.

Now, suppose that $A$, $B$ are disjoint vertex sets in $G$ such that $E_G[A,B]$ is a matching that covers $A \cup B$. Define $D= \overline{A}$ and observe that every $x \in  A$ is super dominated by the vertex $y$ that is the pair of $x$ in the matching $E_G[A,B]$. Indeed, by our assumption, $y \in D$ and $N_G(y)\cap A = \{x\}$. It also follows that $B$ can be considered as a core of $D$. \qed

Note that Lemma~\ref{lem:matching} implies the lower bound of Theorem~\ref{thm:alpha'-rho}. Moreover, it also implies the following. 

\begin{corollary}
\label{cor:exchange}
If $D$ is a super dominating set of a graph $G$, and $D^\ast$ is a core of $D$, then $(D\setminus D^\ast) \cup \overline{D}$ is also a super dominating set of $G$. In particular, if $D$ is a $\gsp$-set, then $(D\setminus D^\ast) \cup \overline{D}$ is also a $\gsp$-set. 
\end{corollary} 

\proof
The first assertion follows directly by Lemma~\ref{lem:matching}. And as $|(D\setminus D^\ast) \cup \overline{D}| = |D|$, the second assertion also follows. 
\qed

\begin{corollary}
\label{cor:each-vertex-in-gsp-set}
If $G$ is a graph and $v\in V(G)$, then there exists a $\gsp$-set of $G$ that contains $v$.
\end{corollary} 

\proof
Let $D$ be an arbitrary $\gsp$-set of $G$. If $v\in D$ there is nothing to prove. Otherwise consider a core $D^\ast$ of $G$.  Then $v\in (D\setminus D^\ast) \cup \overline{D}$ which is a $\gsp$-set by Corollary~\ref{cor:exchange}.
\qed 

It is interesting to compare Corollary~\ref{cor:each-vertex-in-gsp-set} with~\cite[Proposition 2.4]{Zhu} which asserts that if $v$ is a leaf of a non-trivial tree $T$, then there exists a $\gsp$-set of $T$ which does not contain $v$. 

\section{Super domination number of trees}
\label{sec:trees}

The main result of this section reads as follows. 

\begin{theorem} \label{thm:tree}
If $T$ is a tree, then $\gsp(T)= n(T)- \alpha'(T)$. Moreover, a $\gsp$-set of $T$ can be determined in linear time over the class of trees.  
\end{theorem}

\proof
Consider a maximum matching $M$ in $T$. We will show that $V(M)$ can be partitioned into two vertex sets $A$ and $B$ such that $E_T[A,B]=M$. By Lemma~\ref{lem:matching}, it will imply that there is a super dominating set of cardinality $n(T)-|A|= n(T)- |M|= n(T)-\alpha'(T)$. Together with the inequality $\gsp(G) \geq n(G)- \alpha'(G)$ from Theorem~\ref{thm:general-bounds} we obtain $\gsp(T) = n(T)- \alpha'(T)$ for the tree, as stated.

To construct the sets $A$ and $B$, we first specify a root vertex $r$  such that $r$ is covered by $M$. We first put $r$ into $A$ and consider the children $v_1, \dots , v_k$ of $r$. If $rv_i \in M$, we put $v_i$ into $B$; if $rv_i \notin M$ but $v_i \in V(M)$, we put $v_i$ into $A$; if $v_i \notin V(M)$,  then $v_i$ remains outside $A \cup B$. We continue analogously while traversing the tree in preorder. When we decide about the children $u_1, \dots , u_\ell$ of a vertex $u$, we have three main cases.
\begin{itemize}
\item First, let  $u \notin V(M)$. Then, if $u_i \in V(M)$, we put $u_i$ into $A$; if $u_i \notin V(M)$, we put it into neither $A$ nor $B$.
\item Suppose that $u \in A$. If $uu_i \in M$, put $u_i$ into $B$; if $uu_i \notin M$ and $u_i \in V(M)$, put $u_i$ into $A$; if $u_i \notin V(M)$, leave $u_i$ outside $A \cup B$.
\item The case when $u \in B$ is analogous to the previous one. If $uu_i \in M$, we put $u_i$ into $A$; if $uu_i \notin M$ and $u_i \in V(M)$, we put $u_i$ into $B$; if $u_i \notin V(M)$, we leave $u_i$ outside $A \cup B$.
\end{itemize}
It is clear that for the constructed sets, $(A,B)$ results in a partition of $V(M)$ such that $M=E_T[A,B]$. Therefore, by  Lemma~\ref{lem:matching}, $D= \overline{A}$ is a super dominating set in $T$ and we may infer $\gsp(T) = n(T)- \alpha'(T)$.

Concerning the construction of a $\gsp$-set of a tree, we remark that a maximum matching of a tree can be obtained in linear time. Once the matching $M$ is in hand, the algorithm described in the proof assigns labels $A$, $B$, $\overline{A\cup B}$ to the vertices in preorder, visiting every vertex only once and making a choice according to local properties. Thus, the determination of a $\gsp$-set of a tree can be done in linear time as stated. \qed

Extending the definition of the subdivision of a graph by setting $S_0(G) = G$, the following result can be considered as a generalization of Theorem~\ref{thm:paths-cycles-stars}(iii).  

\begin{corollary}
\label{cor:stars}
If $k\ge 0$ and $n\ge 2$, then 
\begin{displaymath}
\gsp(S_k(K_{1,n})) = \left\{ \begin{array}{ll}
\frac{n(k+2)}{2}; & k\ \mbox{even},\\ \\
\frac{n(k+1)}{2} + 1; & k\ \mbox{odd}.
\end{array} \right.
\end{displaymath}
\end{corollary}

\proof
By Theorem~\ref{thm:paths-cycles-stars}(iii), $\gsp(S_0(K_{1,n})) = \gsp(K_{1,n}) = n$, hence the assertion holds for $k=0$. 

It is straightforward to see that if $k\ge 2$ is even, then $\alpha'(S_k(K_{1,n})) = n\frac{k}{2} + 1$, and if $k\ge 1$ is odd, then $\alpha'(S_k(K_{1,n})) = n\frac{k+1}{2}$. The result now follows by applying Theorem~\ref{thm:tree}. 
\qed

\section{Bipartite graphs}
\label{sec:bipartite}

By Theorem~\ref{thm:tree}, the super domination number and a $\gsp$-set can be determined in linear time for trees. In this section we show that the same problem is NP-hard over the class of bipartite graphs.

\begin{theorem}	\label{thm:bipartite}
	\begin{itemize}
		\item[$(a)$] 	It is NP-complete to decide whether $ \gsp(G) \leq k$  holds if $G$ is a bipartite graph of girth $g(G) \ge 8$ and the positive integer $k$  is part of the input.
		\item[$(b)$] It is NP-complete to decide whether $\gsp(G) = n(G)-\alpha'(G)$  holds if $G$ is bipartite and $g(G)\ge 8$. 
	\end{itemize}
\end{theorem}

\proof As $\gsp(G) \geq n(G)-\alpha'(G)$ holds for every graph $G$, the equality in $(b)$ is equivalent to the inequality $\gsp(G) \leq n(G)-\alpha'(G)$. Thus, both decision problems $(a)$ and $(b)$ belong to NP. In order to prove that the decision problems in $(a)$ and $(b)$ are NP-hard, we present a polynomial-time reduction from 3-SAT problem, which is a classical NP-complete problem \cite{Garey1979}.
\medskip

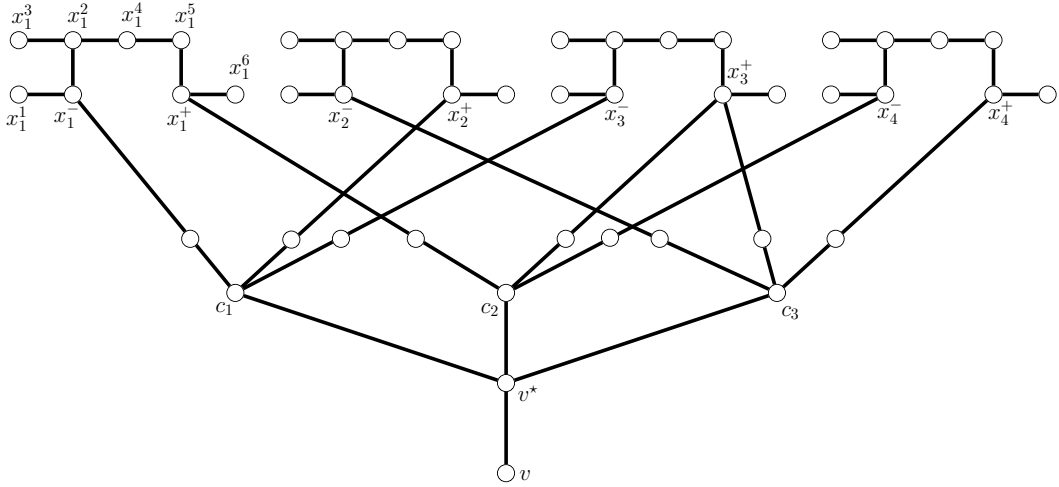
\begin{figure}[ht!] \label{fig:bipartite}
	\begin{center}
		\psscalebox{0.6 0.6}
		{
			\begin{pspicture}(0,-9.475695)(23.28139,1.5756946)
				\psline[linecolor=black, linewidth=0.08](0.28,0.3256946)(1.48,0.3256946)(1.48,-0.8743054)(0.28,-0.8743054)(0.28,-0.8743054)
				\psline[linecolor=black, linewidth=0.08](1.48,0.3256946)(2.68,0.3256946)(3.88,0.3256946)(3.88,-0.8743054)(5.08,-0.8743054)(5.08,-0.8743054)
				\psline[linecolor=black, linewidth=0.08](6.28,0.3256946)(7.48,0.3256946)(7.48,-0.8743054)(6.28,-0.8743054)(6.28,-0.8743054)
				\psline[linecolor=black, linewidth=0.08](7.48,0.3256946)(8.68,0.3256946)(9.88,0.3256946)(9.88,-0.8743054)(11.08,-0.8743054)(11.08,-0.8743054)
				\psline[linecolor=black, linewidth=0.08](12.28,0.3256946)(13.48,0.3256946)(13.48,-0.8743054)(12.28,-0.8743054)(12.28,-0.8743054)
				\psline[linecolor=black, linewidth=0.08](13.48,0.3256946)(14.68,0.3256946)(15.88,0.3256946)(15.88,-0.8743054)(17.08,-0.8743054)(17.08,-0.8743054)
				\psline[linecolor=black, linewidth=0.08](18.28,0.3256946)(19.48,0.3256946)(19.48,-0.8743054)(18.28,-0.8743054)(18.28,-0.8743054)
				\psline[linecolor=black, linewidth=0.08](19.48,0.3256946)(20.68,0.3256946)(21.88,0.3256946)(21.88,-0.8743054)(23.08,-0.8743054)(23.08,-0.8743054)
				\psline[linecolor=black, linewidth=0.08](1.48,-0.8743054)(5.08,-5.2743053)(5.08,-5.2743053)
				\psline[linecolor=black, linewidth=0.08](5.08,-5.2743053)(9.88,-0.8743054)(9.88,-0.8743054)
				\psline[linecolor=black, linewidth=0.08](3.88,-0.8743054)(11.08,-5.2743053)(11.08,-5.2743053)
				\psline[linecolor=black, linewidth=0.08](5.08,-5.2743053)(13.48,-0.8743054)(13.48,-0.8743054)
				\psline[linecolor=black, linewidth=0.08](11.08,-5.2743053)(15.88,-0.8743054)(15.88,-0.8743054)
				\psline[linecolor=black, linewidth=0.08](11.08,-5.2743053)(19.48,-0.8743054)(19.48,-0.8743054)
				\psline[linecolor=black, linewidth=0.08](7.48,-0.8743054)(17.08,-5.2743053)(17.08,-5.2743053)
				\psline[linecolor=black, linewidth=0.08](15.88,-0.8743054)(17.08,-5.2743053)(17.08,-5.2743053)
				\psline[linecolor=black, linewidth=0.08](17.08,-5.2743053)(21.88,-0.8743054)(21.88,-0.8743054)
				\psline[linecolor=black, linewidth=0.08](5.08,-5.2743053)(11.08,-7.2743053)(11.08,-5.2743053)(11.08,-5.2743053)(11.08,-5.2743053)
				\psline[linecolor=black, linewidth=0.08](11.08,-7.2743053)(17.08,-5.2743053)(17.08,-5.2743053)
				\psline[linecolor=black, linewidth=0.08](11.08,-9.274305)(11.08,-7.2743053)(11.08,-7.2743053)
				\psdots[linecolor=black, fillstyle=solid, dotstyle=o, dotsize=0.4, fillcolor=white](11.08,-9.274305)
				\psdots[linecolor=black, fillstyle=solid, dotstyle=o, dotsize=0.4, fillcolor=white](11.08,-7.2743053)
				\psdots[linecolor=black, fillstyle=solid, dotstyle=o, dotsize=0.4, fillcolor=white](11.08,-5.2743053)
				\psdots[linecolor=black, fillstyle=solid, dotstyle=o, dotsize=0.4, fillcolor=white](5.08,-5.2743053)
				\psdots[linecolor=black, fillstyle=solid, dotstyle=o, dotsize=0.4, fillcolor=white](17.08,-5.2743053)
				\psdots[linecolor=black, fillstyle=solid, dotstyle=o, dotsize=0.4, fillcolor=white](4.08,-4.0743055)
				\psdots[linecolor=black, fillstyle=solid, dotstyle=o, dotsize=0.4, fillcolor=white](6.32,-4.0943055)
				\psdots[linecolor=black, fillstyle=solid, dotstyle=o, dotsize=0.4, fillcolor=white](7.42,-4.0743055)
				\psdots[linecolor=black, fillstyle=solid, dotstyle=o, dotsize=0.4, fillcolor=white](9.08,-4.0743055)
				\psdots[linecolor=black, fillstyle=solid, dotstyle=o, dotsize=0.4, fillcolor=white](12.4,-4.0743055)
				\psdots[linecolor=black, fillstyle=solid, dotstyle=o, dotsize=0.4, fillcolor=white](13.38,-4.0543056)
				\psdots[linecolor=black, fillstyle=solid, dotstyle=o, dotsize=0.4, fillcolor=white](14.48,-4.0743055)
				\psdots[linecolor=black, fillstyle=solid, dotstyle=o, dotsize=0.4, fillcolor=white](16.76,-4.0743055)
				\psdots[linecolor=black, fillstyle=solid, dotstyle=o, dotsize=0.4, fillcolor=white](18.38,-4.0743055)
				\psdots[linecolor=black, fillstyle=solid, dotstyle=o, dotsize=0.4, fillcolor=white](0.28,-0.8743054)
				\psdots[linecolor=black, fillstyle=solid, dotstyle=o, dotsize=0.4, fillcolor=white](0.28,0.3256946)
				\psdots[linecolor=black, fillstyle=solid, dotstyle=o, dotsize=0.4, fillcolor=white](1.48,0.3256946)
				\psdots[linecolor=black, fillstyle=solid, dotstyle=o, dotsize=0.4, fillcolor=white](1.48,-0.8743054)
				\psdots[linecolor=black, fillstyle=solid, dotstyle=o, dotsize=0.4, fillcolor=white](2.68,0.3256946)
				\psdots[linecolor=black, fillstyle=solid, dotstyle=o, dotsize=0.4, fillcolor=white](3.88,0.3256946)
				\psdots[linecolor=black, fillstyle=solid, dotstyle=o, dotsize=0.4, fillcolor=white](3.88,-0.8743054)
				\psdots[linecolor=black, fillstyle=solid, dotstyle=o, dotsize=0.4, fillcolor=white](5.08,-0.8743054)
				\psdots[linecolor=black, fillstyle=solid, dotstyle=o, dotsize=0.4, fillcolor=white](6.28,-0.8743054)
				\psdots[linecolor=black, fillstyle=solid, dotstyle=o, dotsize=0.4, fillcolor=white](6.28,0.3256946)
				\psdots[linecolor=black, fillstyle=solid, dotstyle=o, dotsize=0.4, fillcolor=white](7.48,0.3256946)
				\psdots[linecolor=black, fillstyle=solid, dotstyle=o, dotsize=0.4, fillcolor=white](7.48,-0.8743054)
				\psdots[linecolor=black, fillstyle=solid, dotstyle=o, dotsize=0.4, fillcolor=white](9.88,0.3256946)
				\psdots[linecolor=black, fillstyle=solid, dotstyle=o, dotsize=0.4, fillcolor=white](9.88,-0.8743054)
				\psdots[linecolor=black, fillstyle=solid, dotstyle=o, dotsize=0.4, fillcolor=white](11.08,-0.8743054)
				\psdots[linecolor=black, fillstyle=solid, dotstyle=o, dotsize=0.4, fillcolor=white](12.28,-0.8743054)
				\psdots[linecolor=black, fillstyle=solid, dotstyle=o, dotsize=0.4, fillcolor=white](13.48,-0.8743054)
				\psdots[linecolor=black, fillstyle=solid, dotstyle=o, dotsize=0.4, fillcolor=white](13.48,0.3256946)
				\psdots[linecolor=black, fillstyle=solid, dotstyle=o, dotsize=0.4, fillcolor=white](12.28,0.3256946)
				\psdots[linecolor=black, fillstyle=solid, dotstyle=o, dotsize=0.4, fillcolor=white](14.68,0.3256946)
				\psdots[linecolor=black, fillstyle=solid, dotstyle=o, dotsize=0.4, fillcolor=white](15.88,0.3256946)
				\psdots[linecolor=black, fillstyle=solid, dotstyle=o, dotsize=0.4, fillcolor=white](15.88,-0.8743054)
				\psdots[linecolor=black, fillstyle=solid, dotstyle=o, dotsize=0.4, fillcolor=white](17.08,-0.8743054)
				\psdots[linecolor=black, fillstyle=solid, dotstyle=o, dotsize=0.4, fillcolor=white](18.28,-0.8743054)
				\psdots[linecolor=black, fillstyle=solid, dotstyle=o, dotsize=0.4, fillcolor=white](18.28,0.3256946)
				\psdots[linecolor=black, fillstyle=solid, dotstyle=o, dotsize=0.4, fillcolor=white](19.48,0.3256946)
				\psdots[linecolor=black, fillstyle=solid, dotstyle=o, dotsize=0.4, fillcolor=white](19.48,-0.8743054)
				\psdots[linecolor=black, fillstyle=solid, dotstyle=o, dotsize=0.4, fillcolor=white](20.68,0.3256946)
				\psdots[linecolor=black, fillstyle=solid, dotstyle=o, dotsize=0.4, fillcolor=white](21.88,0.3256946)
				\psdots[linecolor=black, fillstyle=solid, dotstyle=o, dotsize=0.4, fillcolor=white](21.88,-0.8743054)
				\psdots[linecolor=black, fillstyle=solid, dotstyle=o, dotsize=0.4, fillcolor=white](23.08,-0.8743054)
				\rput[bl](11.38,-9.434305){\Large{$v$}}
				\rput[bl](11.34,-7.6743054){\Large{$v^\star$}}
				\rput[bl](4.64,-5.7743053){\Large{$c_1$}}
				\rput[bl](10.52,-5.7943053){\Large{$c_2$}}
				\rput[bl](17.18,-5.8543053){\Large{$c_3$}}
				\rput[bl](0.0,-1.6743054){\Large{$x_1^1$}}
				\rput[bl](1.36,0.6056946){\Large{$x_1^2$}}
				\rput[bl](0.12,0.58569455){\Large{$x_1^3$}}
				\rput[bl](2.56,0.6256946){\Large{$x_1^4$}}
				\rput[bl](3.74,0.6056946){\Large{$x_1^5$}}
				\rput[bl](4.94,-0.5343054){\Large{$x_1^6$}}
				\rput[bl](3.56,-1.6943054){\Large{$x_1^+$}}
				\rput[bl](1.06,-1.6143054){\Large{$x_1^-$}}
				\rput[bl](9.78,-1.5943054){\Large{$x_2^+$}}
				\rput[bl](15.98,-0.6343054){\Large{$x_3^+$}}
				\rput[bl](21.76,-1.5543054){\Large{$x_4^+$}}
				\rput[bl](7.14,-1.6143054){\Large{$x_2^-$}}
				\rput[bl](13.26,-1.5743054){\Large{$x_3^-$}}
				\rput[bl](19.3,-1.5343055){\Large{$x_4^-$}}
				\psdots[linecolor=black, dotstyle=o, dotsize=0.4, fillcolor=white](8.68,0.3256946)
			\end{pspicture}
		}
	\end{center}
	\caption{Graph $G_F$ for the formula $F=(\neg x_1 \vee x_2 \vee \neg x_3) \wedge (x_1 \vee x_3 \vee \neg x_4) \wedge (\neg x_2 \vee x_3 \vee x_4) $ constructed in the proof of Theorem~\ref{thm:bipartite} }
	\label{NP-completeness}
\end{figure}

Let $F$ be a $3$-SAT instance with clauses $ C_1,\dots, C_\ell $ over the Boolean variables $ x_1, \dots, x_s$. We construct a graph $G_F$ such that $F$ is satisfiable if and only if $\gsp(G_F)\leq 4s+3\ell +1$. 
\paragraph{Construction of $G_F$.} For each variable $x_i$, we take eight vertices that form the set $X_i= \{x_i^-, x_i^+, x_i^1, \dots , x_i^6\}$ and add edges such that $x_i^1x_i^-x_i^2x_i^4x_i^5x_i^+x_i^6$ is an induced path  and $x_i^3x_i^4$ is a pendant edge in $G_F$.  Each clause $ C_j $, will be represented by a vertex $ c_j $ in $G_F$. If $x_i$ is a literal in $C_j$, we add an edge $x_i^+c_j$  and subdivide it by a vertex $y_{j,i}$. Similarly, if $\neg x_i$ is a literal in $C_j$, we add an edge $x_i^-c_j$  and subdivide it by a vertex $y_{j,i}$.
The set of these subdivision vertices will be denoted by $Y$. To finish the construction, we add two further vertices, namely $v$ and $v^*$, and the edges $vv^*$ and $v^*c_j$ for each $j \in [\ell]$. (See Fig.~\ref{fig:bipartite} for illustration.) It is easy to check that the constructed graph $G_F$ is bipartite with $n(G_F)=8s+4\ell+2$ and, moreover, if $G$ is not a tree\footnote{We may suppose, without loss of generality, that $G_F$ is of girth of at least $8$. Indeed, if $G_F$ is a tree, we may consider the formula $F'=F \wedge (x_1 \vee \neg x_1 \wedge x_2) $. Clearly, a truth assignment satisfies $F$  if and only if it satisfies $F'$, and $\gsp(G_F)= n(G_F)-\alpha'(G_F)$  if and only if the same is true for $G_{F'}$.}, then its girth is at least $8$. 

We first prove that $\alpha'(G_F)=4s+\ell +1$. Let $M$ be a matching in $G_F$. For every $i \in [s]$, the vertices in $X_i$ may be incident with at most four edges from $M$. Each clause vertex $c_j$ may be incident with one edge from $M$. The only edge in $G_F$ that is not covered by the previous vertices is $vv^*$ and it may belong to $M$ only if  $c_jv^* \notin M$ holds for all $j \in [\ell]$. This proves $\alpha'(G_F) \leq 4s + \ell +1$ and it is easy to find a matching of size $4s + \ell +1$ in $G_F$. Therefore, $\alpha'(G_F) = 4s + \ell +1$. Moreover, every maximum matching contains the following edges: $vv^*$; $x_i^2x_i^3$ and $x_i^4x_i^5$ for every $i \in [s]$; one edge between $c_j$ and $Y$ for every $j \in [\ell]$; one edge between $x_i^+$ and $Y \cup \{x_i^6\}$ and one edge between $x_i^-$ and $Y \cup \{x_i^1\}$ for every $i \in [s]$. 

Now, suppose that $\gsp(G_F)\leq 4s+3\ell +1$ holds and prove that the $3$-SAT formula $F$ is satisfiable. Since $\alpha'(G_F)=4s+\ell +1$, by Theorem~\ref{thm:general-bounds} the condition is equivalent with $\gsp(G_F)= 4s+3\ell +1$. Let $D$ be a minimum $\gsp$-set in $G_F$. By Lemma~\ref{lem:matching}, there are two disjoint vertex sets $A=\overline{D}$ and $B=D^\ast$ such that $|A|=|B|= n(G_F)- \gsp(G_F)= 4s+\ell +1$ and $E_{G_F}[A,B]$ is a matching $M$. Since $M$ is a maximum matching, for every $i \in [s]$, we have  $x_i^-, x_i^+ \in V(M)$ and $x_i^2x_i^3,  x_i^4x_i^5 \in M$. By the condition $E_{G_F}[A,B]=M$, if $x_i^+ \in A$ holds, then $x_i^5 \in A$ and $x_i^4, x_i^2, x_i^- \in B$ follow. Analogously, if $x_i^+ \in B$, we may conclude $x_i^- \in A$. Therefore, we may define a truth function $\varphi\colon X \to \{\mbox{true, false}\}$ in the following way: 
\begin{displaymath}
	\varphi(x_i)= \left\{ \begin{array}{ll}
		\mbox{true}; & x_i^+ \in B, \\
		\\
		\mbox{false}; & x_i^- \in B.
	\end{array} \right.
\end{displaymath}

By Corollary~\ref{cor:exchange}, we may suppose that $v^* \in A$. 
Consider a clause vertex $c_j$. As  $v^*c_j \notin M$ and $c_j \in V(M)$, the vertex $c_j$ also belongs to $A$. If $y_{j,i}$ is the vertex from $Y$ such that $c_jy_{j,i} \in M$, then $y_{j,i} \in B$. Suppose first that the other neighbor of $y_{j,i}$ is $x_i^+$ i.e., the clause $C_j$ contains the positive literal $x_i$. As $x_i \in V(M)$ and $y_{j,i}$ is already covered by one matching edge, $y_{j,i}x_i^+ \notin M$ and $x_i^+ \in B$ holds. Then, by definition, we have $\varphi(x_i)=\mbox{true}$ and the positive literal $x_i^+$ satisfies clause $C_j$.  Similarly, if the other neighbor of $y_{j,i}$ is $x_i^-$, then  $C_j$ contains the literal $\neg x_i$. As $E_{G_F}[A,B]=M$, we may infer $x_i^- \in B$. It implies $\varphi(x_i)=\mbox{false}$ and hence, the literal $\neg x_i$ satisfies $C_j$. It is true for all clauses in   $F$ and proves the satisfiability of the formula. 
\medskip

To prove the other direction of the statement, we suppose that $F$ is satisfied by a truth assignment $\phi\colon X \to \{\mbox{true, false}\}$. Let us define
$$
		D   =  Y \cup \{v\}  \cup   \{x_i^1, x_i^3, x_i^5, x_i^+ \colon i \in [s] \mbox{ and } \phi(x_i)=\mbox{true}  \} $$
		$$\cup \, \{x_i^-, x_i^2, x_i^4, x_i^6 \colon i \in [s] \mbox{ and } \phi(x_i)=\mbox{false}  \}.
	$$

It is easy to check that $D$ is a super dominating set and $|D|=4s+3\ell+1$. Indeed, it is enough to consider the following connections:
\begin{itemize}
\item $v$ super dominates $v^*$; 
\item if $\phi(x_i)=\mbox{true}$, then $x_i^1$ super dominates $x_i^-$, $x_i^3$ super dominates $x_i^2$, 
$x_i^5$ super dominates $x_i^4$, and 
$x_i^+$ super dominates $x_i^6$;
\item if $\phi(x_i)=\mbox{false}$, then $x_i^-$ super dominates $x_i^1$, $x_i^2$ super dominates $x_i^3$, 
$x_i^4$ super dominates $x_i^5$, and 
$x_i^6$ super dominates $x_i^+$;
\item if a clause $C_j$ is satisfied by a literal $x_i$ or $\neg x_i$, then the corresponding subdivision vertex $y_{j,i} \in D$ and $N(y_{j,i}) \cap \overline{D}= \{c_j \}$ and thus, $y_{j,i}$ super dominates $c_j$.
\end{itemize}

We have proved that the NP-complete problem 3-SAT can be reduced to the problem of deciding whether $\gsp(G_F) \leq n(G_F)- \alpha'(G_F)= 4s+3\ell +1$ holds. The reduction can be done in polynomial time and therefore, we may conclude that both problems $(a)$ and $(b)$ are NP-complete.   
\qed

\section{Super domination in subdivision graphs}
\label{sec:subdivisions}

\subsection{$(4t+3)$-subdivisions}

For a graph $G$, let $\widehat{n}(G)$ be the maximum size of a subset $\widehat{V} \subseteq V(G)$ such that there exists an injective mapping $\phi: \widehat{V} \to E(G)$ so that $v \in \phi(v)$ holds for every $v \in \widehat{V}$. We will say that a function $\phi$ with these properties is a \emph{DR-function} in $G$; and if $\widehat{n} (G)=n(G)$, we may say that the vertex set of $G$ has a \emph{set of distinct representatives (SDR)}.
\begin{lemma} \label{lem:SDR}
If $G$ is a connected graph that is not a tree, then $\widehat{n}(G)=n(G)$. If $G$ is a tree, then $\widehat{n}(G)=n(G)-1$.
\end{lemma}
\proof Associate every vertex $v \in V(G)$ with the set $E(v)$ of edges that are incident to $v$. First, consider a proper subset $X $ of $V(G)$ and the set $E(X)=\bigcup_{v \in X}  E(v)$. In the induced subgraph $G[X]$, every component $F$ satisfies $m(F) \ge n(F)-1$ and, since $G$ is connected and $F \neq G$, the vertex set of $F$ is incident with at least one edge not contained in the subgraph $F$. These extra edges are pairwise different for different components of $G[X]$. We therefore conclude $|E(X)|\ge |X|$ for every vertex set $X \subsetneqq V(G)$. Notice that it is true for every graph $G$, no matter $G$ is a tree or not. Consider now the case $X=V(G)$. If $G$ is not a tree, then $n(G)=|X| \leq |E(X)|=m(G)$ and, as Hall's Condition is satisfied for each $X \subseteq V(G)$, there exists a system of distinct representatives for the vertex set of $G$. That is, $\widehat{n}(G)=n(G)$. If $G$ is a tree, $|V(G)| > |E(G)|$ and there is no SDR for the vertex set. On the other hand, if we consider $G$ as a tree rooted in $r$ and map every non-root vertex $v$ to the edge between $v$ and its parent, the obtained mapping is a DR-function from $V(G)\setminus \{r\}$ to $E(G)$. It proves $\widehat{n}(G)= n(G)-1$ for every tree $G$. \qed

\begin{theorem}
	\label{prop:4t+3}
	For every connected graph $G$ and integer $k \equiv 3 \bmod 4$, 
	\begin{displaymath}
	\gsp(S_k(G))= \left\{ \begin{array}{ll}
	\frac{k+1}{2}\, m(G)+1; & \mbox{$G$ is a tree}, \\
	\\
	\frac{k+1}{2}\, m(G); & \mbox{otherwise}.
	\end{array} \right.
	\end{displaymath}
\end{theorem}
\proof
Let $V(G)=\{v_1,\ldots,v_n\}$ and $k=4t+3$. Observe that $n(S_k(G))=n+k\,m(G)$.
We first show that 
\begin{eqnarray} \label{ineq:alpha'}
	\alpha'(S_k(G)) \leq \frac{k-1}{2}\, m(G) +\widehat{n}(G).
\end{eqnarray}

Let $M$ be a maximum matching in $S_k(G)$. For every super edge $P_{v_iv_j}$, we have two possibilities:
\begin{itemize}
\item[$(a)$] $M$ contains at most $\frac{k-1}{2}$ edges from $P_{v_iv_j}$;
\item[$(b)$] $M$ contains exactly $\frac{k+1}{2}$ edges from $P_{v_iv_j}$ and at least one of $v_i$ and $v_j$ is covered by a matching edge belonging to $P_{v_iv_j}$. 
\end{itemize}
As each $v_i \in V(G)$ is covered by at most one edge from $M$, the number of super edges satisfying $(b)$ is at most $n$. Moreover, if $(b)$ is valid for a super edge $P_{v_iv_j}$, then $v_i(v_iv_j)_1 \in M$ or $(v_iv_j)_k v_j \in M$. In the first case, we define  $\phi(v_i) =v_iv_j$, while we set $\phi(v_j) =v_iv_j$ in the latter case. (If both edges $v_i(v_iv_j)_1$ and  $(v_iv_j)_k v_j $ belong to $M$, then to keep $\phi$ injective, we set just $\phi(v_i) =v_iv_j$ for the smaller index $i$.) As $\phi$ is a DR-function, the number of super edges with property $(b)$ is at most $\widehat{n}(G)$. This proves the inequality (\ref{ineq:alpha'}), and together with Theorem~\ref{thm:alpha'-rho} we conclude 
\begin{eqnarray} \label{ineq:2alpha'}
\gsp(S_k(G)) \geq n(S_k(G))- \alpha'(S_k(G)) 
\geq \frac{k+1}{2}\, m(G) + (n(G)- \widehat{n}(G)) 
\end{eqnarray}
where, according to Lemma~\ref{lem:SDR}, the last term is $1$ if $G$ is a tree and $0$ if $G$ contains a cycle. 
\medskip

To prove the other direction, we construct a $\gsp$-set $D$  for $G$.
Let $\phi$ be a DR-function of $G$ with domain $V(G)$ if $G$ is not a tree, and with domain $V(G)\setminus \{v_n\}$ otherwise.
\begin{itemize}
\item[$(i)$] If $\phi(v_i)=v_iv_j$, let $D$ contain the following vertices from the super edge $P_{v_iv_j}$:
$$ (v_iv_j)_1, (v_iv_j)_2, (v_iv_j)_5, (v_iv_j)_6, \dots, (v_iv_j)_{4t+1},
(v_iv_j)_{4t+2}.
$$
\item[$(ii)$] If $v_iv_j$ does not belong to the image set of $\phi$ and $i <j$, let $D$ contain the following vertices from the super edge $P_{v_iv_j}$:
$$ (v_iv_j)_2, (v_iv_j)_3, (v_iv_j)_6, (v_iv_j)_7, \dots, (v_iv_j)_{4t+2},
(v_iv_j)_{4t+3}.
$$
\item[$(iii)$] If $G$ is a tree and $v_n$ does not have a representative edge in $\phi$, then $v_n$ also belongs to $D$. Note that the other vertices of $G$ belong to $\overline{D}$.
\end{itemize}

If an internal subdivision vertex, $(v_iv_j)_s$ with $3 \leq s \leq 4t+1$, does not belong to $D$, it is easy to identify a neighbor that super dominates it. A vertex $v_i$ with $\phi(v_i)=v_iv_j$ is always super dominated by $(v_iv_j)_1$. The subdivision vertices $(v_iv_j)_2$ and $(v_iv_j)_{4t+2}$ always belong to $D$; if $(v_iv_j)_1 \notin D$, it is super dominated by $(v_iv_j)_2$; if $(v_iv_j)_{4t+3} \notin D$, it is super dominated by $(v_iv_j)_{4t+2}$.

 No matter whether $(i)$ or $(ii)$ was applied when we specified the vertices in $V(P_{v_iv_j}) \cap D$, we added exactly $2t+2= \frac{k+1}{2}$ subdivision vertices to $D$ in each step. Thus, $D$ contains $\frac{k+1}{2}\, m(G)$ subdivision vertices and also contains $v_n$ if $G$ is a tree. This proves the upper bound
$$\gsp(S_k(G)) \leq 
 \frac{k+1}{2}\, m(G) + (n(G)- \widehat{n}(G)).
 $$
  We infer that the equality $\gsp(S_k(G)) =
  \frac{k+1}{2}\, m(G) + (n(G)- \widehat{n}(G))$ holds for every graph $G$ as stated. \qed

For connected graphs, Theorem~\ref{prop:4t+3} and inequality (\ref{ineq:2alpha'}) in its proof together imply the following statement. Since $\gsp(F)$ and $\alpha'(F)$ are additive under disjoint union of graphs, we may state:
\begin{proposition}
For every graph $G$ and integer $k \equiv 3 \bmod 4$, it holds that 	
$$\gsp(S_k(G)) = n(S_k(G))- \alpha'(S_k(G)).$$ 
\end{proposition}
As the number of tree components in $G$ can be computed in linear time, and $\gsp(F)$ is additive under taking disjoint union of graphs, we conclude the subsection with the following consequence of Theorem~\ref{prop:4t+3}.

\begin{theorem} \label{comp:4t+3}
If $k$ is a positive integer with $k \equiv 3 \bmod 4$, then the super domination number   can be computed in linear time over the class of $k$-subdivision graphs.
\end{theorem}

\subsection{$(4t+1)$-subdivisions}

\begin{theorem}
	\label{prop:4t+1}
	For every connected graph $G$ and integer $k \equiv 1 \bmod 4$, 
	\begin{displaymath}
	\gsp(S_k(G))= \left\{ \begin{array}{ll}
	\frac{k+1}{2}\, m(G); & \mbox{$G$ contains an even cycle}; \\
	\\
	\frac{k+1}{2}\, m(G)+1; & \mbox{otherwise}.
	\end{array} \right.
	\end{displaymath}
\end{theorem}
\proof
Let $V(G)=\{v_1,\ldots,v_n\}$ and $k=4t+1$. Suppose first that $D$ is a minimum super dominating set of $S_k(G)$ and consider $A=\overline{D}$ and a core $B$ of $D$. By Lemma~\ref{lem:matching}, the edges $E[A,B]$ form a matching $M$ in $S_k(G)$. If $M$ is fixed, we have three possibilities for a super edge $P_{v_iv_j}$.
\begin{itemize}
	\item[$(a)$] $M$ contains at most $\frac{k-1}{2}$ edges from $P_{v_iv_j}$. The set of the corresponding edges $v_iv_j \in E(G)$ will be denoted by $E_0$. 
	\item[$(b)$] $M$ contains exactly $\frac{k+1}{2}$ edges from $P_{v_iv_j}$ and exactly one of $v_i$ and $v_j$ is covered by an edge from $M \cap E(P_{v_iv_j})$. If this vertex, say $v_i$, is contained in $A$, we set $v_iv_j \in E_A$. Similarly, if $v_i \in B$, $v_i(v_iv_j)_1 \in M$, and $(v_iv_j)_kv_j \notin M$, then the edge $v_iv_j$ belongs to $E_B$.
	\item[$(c)$] $M$ contains exactly $\frac{k+1}{2}$ edges from $P_{v_iv_j}$ and both $v_i(v_iv_j)_1$ and $(v_iv_j)_kv_j$ belong to $M$. In this case, we set $v_iv_j \in E_{2}$.	 
\end{itemize}
As $E_0, E_A, E_B, E_{2}$ is a partition of $E(G)$, we may estimate the size of $M$ as follows:
\begin{eqnarray} \label{ineq1:4t+1}
|M|\leq |E_0|\,  \frac{k-1}{2}+(|E_A|+|E_B|+|E_2|)\,  \frac{k+1}{2}= \frac{k-1}{2}\, m(G) +|E_A|+|E_B|+|E_2|.
\end{eqnarray}
By definition, if $v_iv_j \in E_A \cup E_B$, then only one of $v_i(v_iv_j)_1$ and $(v_iv_j)_kv_j$ belongs to $M$. If $v_iv_j \in E_2$, then both $v_i(v_iv_j)_1$ and $(v_iv_j)_kv_j$ are contained in $M$. Since each vertex $v_i \in V(G)$ is covered by at most one $M$-edge, we infer $|E_A|+|E_B|+2 |E_2| \leq n$ and, in turn, we get from (\ref{ineq1:4t+1}) that
\begin{eqnarray} \label{ineq2:4t+1}
|M|\leq \frac{k-1}{2}\, m(G) +n - |E_2|.
\end{eqnarray} 
As $|M|=|A|= |\overline{D}|= n(S_k(G))-\gsp(S_k(G))$ and $n(S_k(G))=n+k \,m(G)$,  inequality (\ref{ineq2:4t+1}) implies
\begin{eqnarray} \label{ineq3:4t+1}
\gsp(S_k(G)) =  n(S_k(G))- |M| \geq 
\frac{k+1}{2}\, m(G) + |E_2|. 
\end{eqnarray}
If $E_2\neq\emptyset$ or $G$ contains an even cycle, (\ref{ineq3:4t+1}) itself proves the required lower bound. From now on, we assume that there is no even cycle in $G$ and that $E_2=\emptyset$.

Consider an edge $v_iv_j \in E_A$ with $v_i(v_iv_j)_1 \in M$. As $M$ contains $\frac{k+1}{2}$ edges from  $P_{v_i,v_j}$ that includes $v_i(v_iv_j)_1$ but not $(v_iv_j)_kv_j$, $M$ contains the following edges from the super edge:
$$v_i(v_iv_j)_1, (v_iv_j)_2(v_iv_j)_3, \dots , (v_iv_j)_{4t}(v_iv_j)_{4t+1}.$$
By Lemma~\ref{lem:matching}, $v_i \in A $ implies $(v_iv_j)_{1} \in B$;
the latter implies $(v_iv_j)_{2} \in B$. Since  $(v_iv_j)_{2}(v_iv_j)_{3} \in M$, we infer $(v_iv_j)_{3} \in A$; and so on. We obtain that $(v_iv_j)_s \in B$ if and only if $s \bmod 4 \in \{1,2\}$; otherwise, $(v_iv_j)_s \in A$. In the last step, $v_j \in B$ also follows. It can be proved analogously that $v_iv_j \in E_B$ and $v_i \in B$ implies $v_j \in A$. 

Let 
$$V'= \{v_i: \exists j \in [n] \mbox{ s.t. } v_i(v_iv_j)_1 \in M\}, \quad  V''= \{v_j: \exists i \in [n] \mbox{ s.t. } v_i(v_iv_j)_1 \in M\},
$$ 
and define a DR-function $\phi: V' \to E_A \cup E_B$ such that $\phi(v_i)= v_iv_j$ if $v_i(v_iv_j)_1 \in M$. Since $E_2 = \emptyset$, it is an injective function. Moreover, if the edge $v_iv_j$ is in the image of $\phi$, then one from $v_i$ and $v_j$ belongs to $A$ and the other one to $B$. Thus, $\phi$ remains a DR-function on $V'$, if we consider the following bipartite subgraph $F$ instead of $G$. We first take the induced subgraph $G[V_F]$,  where $V_F= V' \cup V''$, and then delete the edges inside $V_F\cap A$ and $V_F\cap B$. We may also say that this graph $F$ is defined by the edge set $E_G[V_F\cap A, V_F\cap B]$.
 By supposition, $G$ contains no even cycle. Therefore, the bipartite subgraph $F$ contains no cycle at all. By Lemma~\ref{lem:SDR}, $\widehat{n}(F) \leq n(F)-1 \leq n-1$. Consequently, no more than $n-1$ vertices of $G$ are covered by an edge from $M$ in $S_k(G)$. This implies $|E_A| + |E_B| \leq n-1$ and, by (\ref{ineq1:4t+1}), we infer $|M| \leq \frac{k-1}{2}\, m(G) +n- (n-1)$ that yields
$$\gsp(S_k(G)) \geq \frac{k+1}{2}\, m(G) + 1, $$
if $G$ contains no even cycle.
\medskip

To prove the reverse inequalities, we construct a $\gsp$-set $D$  for $G$. First we define an appropriate DR-function $\phi $. 
\begin{itemize}
	\item If $G$ contains an even cycle $C$, take a unicyclic spanning subgraph $H$ of $G$ such that $C$ is the only cycle in $H$. Then $H$ is bipartite, not a tree, and therefore, by Lemma~\ref{lem:SDR}, $\widehat{n}(H)=n$. Let $A_H$ and $B_H$ be the partite classes of $H$. By Lemma~\ref{lem:SDR}, there is a DR-function $\phi$ which assigns a representative edge from $E(H)$ to each vertex from $V(G)$.
	\item If $G$ contains no even cycle, choose a spanning tree $H$ in $G$. Again, $H$ is bipartite, but now we have $\widehat{n}(H)=n-1$. Let $A_H$ and $B_H$ be the partite classes of $H$. By Lemma~\ref{lem:SDR}, we can define a DR-function $\phi$ that assigns a representative edge from $E(H)$ to each vertex from $V(G)\setminus\{v_n\}$.
\end{itemize}
Having a DR-function $\phi$ in hand, we define a super dominating set $D$ in $S_k(G)$ with a size that matches the required upper bound.
\begin{itemize}
	\item[$(i)$] If $v_i \in B_H$ or $v_i \notin A_H \cup B_H$, we set $v_i \in D$.
	\item[$(ii)$] If $\phi(v_i)=v_iv_j$ and $v_i \in A_H$, let $D$ contain the following vertices from the super edge $P_{v_iv_j}$:
	$$ (v_iv_j)_1, (v_iv_j)_2, (v_iv_j)_5, (v_iv_j)_6, \dots, (v_iv_j)_{4t-3}, (v_iv_j)_{4t-2}, (v_iv_j)_{4t+1}.
	$$
	\item[$(iii)$] If $\phi(v_i)=v_iv_j$ and $v_i \in B_H$, let $D$ contain the following vertices from  $P_{v_iv_j}$:
	$$ (v_iv_j)_3, (v_iv_j)_4, (v_iv_j)_7, (v_iv_j)_8, \dots, (v_iv_j)_{4t-1}, (v_iv_j)_{4t}.
	$$
	\item[$(iv)$] If $v_iv_j$ does not belong to the image set of $\phi$, $v_i \in A_H$, and $i <j$, then let $D$ contain the following vertices from  $P_{v_iv_j}$:
	$$ (v_iv_j)_2, (v_iv_j)_3, (v_iv_j)_6, (v_iv_j)_7, \dots, (v_iv_j)_{4t-2},
	(v_iv_j)_{4t-1},
	(v_iv_j)_{4t+1}.
	$$
	\item[$(v)$] If $v_iv_j$ does not belong to the image set of $\phi$, $v_i \in B_H$, and $i <j$, let $D$ contain the following vertices from  $P_{v_iv_j}$:
	$$ (v_iv_j)_1, (v_iv_j)_4, (v_iv_j)_5, (v_iv_j)_8, (v_iv_j)_9, \dots, (v_iv_j)_{4t},
		(v_iv_j)_{4t+1}.
	$$
\end{itemize}

If an internal subdivision vertex, $(v_iv_j)_s$ with $3 \leq s \leq 4t-1$, does not belong to $D$, it is easy to see that a neighbor super dominates it. A vertex $v_i \notin D$ with $\phi(v_i)=v_iv_j$ is always super dominated by $(v_iv_j)_1$.  The subdivision vertex  $u=(v_iv_j)_1$ does not belong to $D$, if $(iii)$ or $(iv)$ was applied. In the latter case, $(v_iv_j)_2$ super dominates $u$. In the first case, $v_i \in B_H $ and $(N[v_i]\setminus \{u\}) \subseteq D$ holds by the determination of $D$. The subdivision vertex  $w=(v_iv_j)_2$ is missing from $D$, if $(iii)$ or $(v)$ was applied. In the first case, $(v_iv_j)_3$ super dominates $w$. For the second case, the  condition in $(v)$ ensures that $v_i \in D$. Hence, $w$ is the only neighbor of $(v_iv_j)_1$ which is outside $D$. In case $(ii)$, vertex $(v_iv_j)_{4t}$ is super dominated by $(v_iv_j)_{4t+1}$ as in this case $v_i \in A_H$ implies $v_j \in B_H$ and therefore, we have $v_j \in D$.  In case $(iv)$, vertex $(v_iv_j)_{4t}$ is super dominated by $(v_iv_j)_{4t-1}.$ A vertex $(v_iv_j)_{4t+1}$ is outside $D$, only if $(iii)$ was applied. In this case, $(v_iv_j)_{4t}$ super dominates it.

Finally, we determine the size of $D$. Case $(i)$ puts $n-|A_H|$ vertices into $D$. When the subdivision vertices are considered, we put $\frac{k+1}{2}$ vertices from each $P_{i,j}$ into $D$, except when case $(iii)$ is applied. There we deal with $|B_H|$ super edges putting $\frac{k-1}{2}$ internal vertices into $D$ from each. This gives
$$\gsp(S_k(G))\leq |D| = n-|A_H| + m(G)\, \frac{k+1}{2} -|B_H|.$$
By the determination of the DR-function $\phi$, $|A_H|+ |B_H| = n$ if $G$ contains an even cycle; and $|A_H|+ |B_H| = n-1$ if every cycle in $G$ is of odd order. Substituting these values in the inequality, we get the required upper bounds on $\gsp(S_k(G))$. This finishes the proof of the theorem.
  \qed
  
Let $\mbox{oc}(G)$ denote the number of components in $G$ that contain no even cycles.  Then Theorem~\ref{prop:4t+1} directly implies:
\begin{proposition}
	For every graph $G$ and integer $k \equiv 1 \bmod 4$, it holds that 	
	$$\gsp(S_k(G)) = \frac{k+1}{2} m(G)+\mbox{oc}(G).$$ 
\end{proposition}

A shortest even cycle in a graph can be found in polynomial (actually quadratic) time~\cite{yuster-1997}, hence the number of even-cycle-free components is easy to determine. Thus we may deduce the following result:
\begin{theorem} \label{comp:4t+1}
If $k$ is a positive integer with $k \equiv 1 \bmod 4$, then the super domination number  can be computed  in polynomial time over the class of $k$-subdivision graphs.
\end{theorem}
\subsection{$4t$-subdivisions} \label{sec:4t}

\begin{theorem}
	\label{thm:4t}
	For every graph $G$ and integer $k \equiv 0 \bmod 4$, 
	\begin{displaymath}
	\gsp(S_k(G))= 
	\frac{k}{2}\, m(G)+ \gsp(G).
	\end{displaymath}
\end{theorem}
\proof
Let $V(G)=\{v_1,\ldots,v_n\}$ and $k=4t$. Suppose that $D$ is a $\gsp$-set of $S_k(G)$ and consider $A=\overline{D}$ and a core $B$ of $D$. By Lemma~\ref{lem:matching}, the edges $E[A,B]$ form a matching $M$ in $S_k(G)$.  We have three possibilities for a super edge $P_{v_iv_j}$.
\begin{itemize}
	\item[$(a)$] $M$ contains at most $\frac{k}{2}-1$ edges from $P_{v_iv_j}$. Let $E_0$ denote the set of the edges $v_iv_j \in E(G)$ with this property.
	\item[$(b)$] $M$ contains exactly $\frac{k}{2}$ edges from $P_{v_iv_j}$. The set of the corresponding  edges $v_iv_j$ in $G$ is denoted by $E_1$.
	\item[$(c)$] $M$ contains exactly $\frac{k}{2}+1$ edges from $P_{v_iv_j}$. In this case, both $v_i(v_iv_j)_1$ and $(v_iv_j)_kv_j$ belong to $M$, and we set $v_iv_j \in E_{2}$.	 
\end{itemize}
By definitions given for $E_0, E_1,  E_{2}$, the following is true: 
\begin{eqnarray} \label{ineq1:4t}
|M|\leq \frac{k}{2}\, m(G) + |E_2|-|E_0|. 
\end{eqnarray}

To prove that $|E_2|-|E_0| \leq n-\gsp(G)$, we first consider a super edge $P_{v_iv_j}$ so that $v_iv_j \in E_2$. The edges $$v_i(v_iv_j)_1, (v_iv_j)_2(v_iv_j)_3, \dots , (v_iv_j)_{4t-2}(v_iv_j)_{4t-1}, (v_iv_j)_{4t}v_j$$ are all included in $M$. If $v_i \in A$ then, by Lemma~\ref{lem:matching}, $E[A,B]$ is a matching and the vertices 
$$(v_iv_j)_{1}, (v_iv_j)_{2}, (v_iv_j)_{5}, (v_iv_j)_{6}, \dots, (v_iv_j)_{4t-2}, v_j$$
 are from $B$; while the remaining subdivision vertices belong to $A$. Therefore, $v_i \in A$ implies $v_j \in B$ and, similarly, $v_i \in B$ implies $v_j \in A$. 

If $v_iv_j \in E_2$, $v_pv_q \in E_2$ with $v_i \in A$, $v_q \in B$, and $G$ contains an edge $v_iv_q$, then we say that $v_iv_q$ is a \emph{critical edge}. As $v_i$ and $v_q$ are already covered by $M$-edges, $v_iv_q \notin E_2$. If $v_iv_q \in E_1$, then the $k/2$ edges in $M \cap E(P_{v_iv_q})$ have to be $$(v_iv_q)_1(v_iv_q)_2, (v_iv_q)_3(v_iv_q)_4, \dots , (v_iv_q)_{k-1}(v_iv_q)_k. $$ 
By Lemma~\ref{lem:matching}, $v_i \in A$ implies $(v_iv_q)_{1} \in A$, $(v_iv_q)_{2}, (v_iv_q)_{3} \in B$, $(v_iv_q)_{4}, (v_iv_q)_{5} \in A, \dots, (v_iv_q)_{k} \in A$. Finally, we infer $v_q \in A$ that contradicts the assumption $v_q \in B$. It implies that $v_iv_q \in E_0$. By symmetry, the same is true if $v_i \in B$ and $v_q \in A$ and therefore, every critical edge belongs to $E_0$.

We now prove that the maximum for $|E_2|-|E_0|$  can be attained without the presence of critical edges. Indeed, if $v_iv_q$ is a critical edge such that $v_iv_j, v_pv_q \in E_2$ and $v_i \in A$, $v_q \in B$, then we may perform the following changes in $M$:
\begin{itemize}
\item Remove the edges $M\cap E(P_{v_i,v_j})$ from $M$ and replace them by the complement edge set $E(P_{v_i,v_j}) \setminus M$. By this change, $v_iv_j$ is moved to $E_1$, and $v_i, v_j$ become uncovered by $M$. So, this step decreases $|E_2|$ by $1$.
After this change $v_i \notin A \cup B$ and we can replace the (at most) $\frac{k}{2} -1$ $M$-edges on $P_{v_i,v_q}$ with the following $\frac{k}{2}$ edges: 
$$(v_iv_q)_{1}(v_iv_q)_{2}, (v_iv_q)_{3}(v_iv_q)_{4}, \dots , (v_iv_q)_{4t-1}(v_iv_q)_{4t}$$
such that we put $(v_iv_q)_{1}, (v_iv_q)_{4}, (v_iv_q)_{5}, \dots , (v_iv_q)_{4t}$ into $B$ and the remaining subdivision vertices into $A$. Since $v_q \in B$, this step keeps the property $E[A,B]=M$ and, by Lemma~\ref{lem:matching}, $D= \overline{A}$ is a super dominating set in $S_k(G)$. Note that this modification removes $v_iv_q$ from $E_0$.
\end{itemize}
After applying the described changes, $|E_2|-|E_0|$ remains the same and we have less critical edges than before. Thus, performing the steps iteratively while there is a critical edge, we obtain a matching $M'$ and sets $A', B'$ without critical edges  such that $|E_2|-|E_0|$ remains unchanged. As there are no critical edges, the $E_2$-edges now form a matching $M^*$ in $G$ such that $E[A^*, B^*]=M^*$ for the sets $A^*=A' \cap V(G)$ and $B^*=B'  \cap V(G)$. Applying Lemma~\ref{lem:matching} again, we conclude that $D^*=\overline{A^*}$ is a super dominating set in $G$. This yields
$$|E_2|= |M^*| =|A^*|=n- |D^*| \leq n-\gsp(G).$$
Now inequality chain in (\ref{ineq1:4t}) can be continued and we obtain
$$ |M|\leq \frac{k}{2}\, m(G) + |E_2|-|E_0| \leq \frac{k}{2}\, m(G) +n- \gsp(G)
$$
which, in turn, proves 
\begin{align*}
\gsp(S_k(G))  & =  n(S_k(G)) - |M|\\
              & \geq   \left(n+k \, m(G)\right) - \left(\frac{k}{2}\, m(G) +n- \gsp(G)\right)\\
              & =  \frac{k}{2}\, m(G) + \gsp(G).
\end{align*}

\medskip

In the second part of the proof we construct a super dominating set $D$ of size $\frac{k}{2}\, m(G) + \gsp(G)$ in $S_k(G)$. Let $D^*$ be a $\gsp$-set in $G$ with the corresponding sets $A^*, B^*$ and matching $M^*=E_G[A^*, B^*]$. 
\begin{itemize}
	\item[$(i)$] If $v_i \in D^*$, we set $v_i \in D$.
	\item[$(ii)$] If $v_iv_j \in M^*$ with $v_i \in A^*$ and $v_j \in B^*$, let $D$ contain the following vertices from the super edge $P_{v_iv_j}$:
	$$ (v_iv_j)_1, (v_iv_j)_2, (v_iv_j)_5, (v_iv_j)_6, \dots, (v_iv_j)_{4t-3}, (v_iv_j)_{4t-2}.
	$$
	\item[$(iii)$] If $v_iv_j \notin M^*$ and $v_i \in A^*$ then, as $E_G[A^*, B^*]=M^*$, we have  $v_j \in V(G) \setminus B^*$. Let us put into $D$ the following subdivision vertices  from $P_{v_iv_j}$:
	$$ (v_iv_j)_2, (v_iv_j)_3, (v_iv_j)_6, (v_iv_j)_7, \dots, (v_iv_j)_{4t-2}, (v_iv_j)_{4t-1}.
	$$
	\item[$(iv)$] If $v_iv_j \notin M^*$ and $v_i \in B^*$ hold and also if both $v_i$ and $v_j$ are outside $A^* \cup B^*$, we put into $D$ the following subdivision vertices  from $P_{v_iv_j}$:
	$$ (v_iv_j)_1, (v_iv_j)_4, (v_iv_j)_5,  \dots, (v_iv_j)_{4t-4}, (v_iv_j)_{4t-3}, (v_iv_j)_{4t}.
	$$
	\end{itemize}
 In step $(i)$, we put $|D^*|= \gsp(G)$ non-subdivision vertices into $D$. Then, for each super edge considered in $(ii)-(iv)$, we put exactly $k/2$ subdivision vertices into $D$. As there are no edge $v_iv_j$ in $G$ with $v_i \in A^*$, $v_j \in B^*$ and $v_iv_j \notin M^*$, we treated each super edge of $S_k(G)$ in the steps $(ii)-(iv)$. These sum up $|D|=\frac{k}{2}\, m(G) + \gsp(G)$.
 
 To check that $D$ is a super dominating set is mainly automatic. We note that if $v_i \notin D$, then $v_i \in A^*$ and there is an edge $v_iv_j \in M^*$ which is considered in $(ii)$. Then,  $(v_iv_j)_1$ super dominates $v_i$. We also remark that in step $(ii)$, the vertex $(v_iv_j)_{4t}$ is super dominated by $v_j$ as all the other super edges $P_{v_jv_p}$ being incident to $v_j$ were considered in step $(iv)$. There, subdivision neighbors $(v_jv_p)_1=(v_pv_j)_{4t}$ were put into $D$. For a super edge $P_{v_iv_j}$ that was treated in $(iv)$, the condition implies $v_i, v_j \in D$. Therefore, $(v_iv_j)_{1}$ and $(v_iv_j)_{4t}$ super dominate $(v_iv_j)_{2}$ and $(v_iv_j)_{4t-1}$, respectively. It shows $\gsp(S_k(G)) \leq |D|  =\frac{k}{2}\, m(G) + \gsp(G)$, and together with the first part of the proof give the equality $\gsp(S_k(G))=\frac{k}{2}\, m(G) + \gsp(G)$ as required. \qed 
\medskip
 
The problem of deciding whether $\gsp(F) \leq \ell$ holds, clearly belongs to NP. Let $\ell$ be part of the input of the problem and $k$ be a fixed integer with $k \equiv 0 \bmod 4$. By Theorem~\ref{thm:bipartite}, it is NP-hard to decide whether $\gsp(G) \leq \ell$ holds over the class of all graphs. By Theorem~\ref{thm:4t}, $\gsp(G) \leq \ell$ holds if and only if $\gsp(S_k(G)) \leq \frac{k}{2}\, m(G) + \ell$.
Thus, we may conclude the following:
\begin{theorem} \label{comp:4t}
Over the class of $k$-subdivision graphs, it is NP-complete to decide whether $\gsp(F) \leq \ell$ holds, if $\ell$ is part of the input and  $k$ is a fixed integer with  $k \equiv 0 \bmod 4$. 
\end{theorem}
 
 \subsection{II-matchings} \label{sec:II}
 Before continuing our study with the last case for subdivision graphs, we introduce a graph invariant and prove an additional complexity result.
 
 In a graph $G$, an \emph{induced matching} is a matching $M \subseteq E(G)$ such that the induced subgraph $G[V(M)]$ contains only the edges from $M$. We denote by $i(G)$ the maximum size of an induced matching in $G$. Induced matchings are applicable in network flow problems, secure communication, VLSI design, and elsewhere, cf.~\cite{golumbic-2000, saeedi-2021}. To decide whether $i(G) \ge \ell$ holds is known to be NP-hard in many classes of graphs, say in planar bipartite graphs~\cite{stockmeyer-1982} and in claw-free graphs~\cite{kobler-2003}. For exact algorithms for maximum induced matchings see~\cite{nguyen-2021, xiao-2017}, and for the complexity aspects of the maximum-weight induced matchings and dominating induced matchings see~\cite{klemz-2022, brandstaedt-2022}, respectively, and references therein. 
 
We further say that a matching $M$ is an \emph{II-matching} if $M$ can be partitioned into two induced matchings $M_1$ and $M_2$. The \emph{II-matching number} $\mbox{ii}(G)$ of $G$ is the maximum size of an II-matching in $G$. We prove that the II-matching number is hard to determine.
 
 \begin{proposition} \label{prop:ii}
 It is NP-complete to decide whether $\mbox{ii}(G) \geq \ell$ holds, if $\ell$ is part of the input.
 \end{proposition}
 
\proof
To decide whether the independence number $\alpha(F)$ of a graph $F$ is at least $k$ is a classical NP-complete problem~\cite{karp}. We show a polynomial-time reduction from the decision problem of $\alpha(F)\geq k$ to the problem of $\mbox{ii}(G) \geq 2k$.

\paragraph{Construction.} For every graph $F$, let $G_F$ be the graph constructed on the vertex set $V(G)\times V(K_4)$, where $V(K_4) = [4]$, by making two different vertices $(x,i)$ and $(y,j)$ adjacent in $G_F$ if either $xy \in E(F)$ or $x=y$. (We note in passing that $G_F$ is isomorphic to the lexicographic product $F\circ K_4$.) Let $V(F)=\{v_1, \dots , v_n\}$ and let $V_i$ denote the vertex set $\{v_i\} \times [4]$ in $G_F$.

\paragraph{Reduction.} We show that $\mbox{ii}(G_F)= 2\alpha(F)$ holds for every graph $F$ and therefore, deciding whether $\alpha(F)\geq k$ is equivalent to the problem of whether $\mbox{ii}(G_F)\geq 2k$ holds.

First, consider a maximum independent set $S$ in $F$ and define the edge sets
$$M_1=\{(x,1)(x,2): x\in S\} \quad \mbox{and} \quad M_2=\{(x,3)(x,4): x\in S\}.$$
By definition, $|M_1|=|M_2|=\alpha(F)$, $M_1 \cup M_2$ is a matching in $G_F$ and, since $S$ is an independent set in $F$, both $M_1$ and $M_2$ are induced matchings. It follows that $\mbox{ii}(G_F) \geq 2\alpha(F)$. 

Assume now, that we have a maximum II-matching $M=M_1 \cup M_2$ in $G_F$. If the induced matching $M_1$ contains a \emph{cross edge} $xy$, that is an edge $xy$ with $x \in V_i$, $y \in V_j$ such that $i \neq j$, then  $M_1$ cannot cover any other vertex from the neighborhood of $x$. Equivalently, if $v_iv_{i'}$ is an edge in $F$, then $V_{i'} \cap V(M_1)$ cannot contain a vertex different from $y$. Thus, if the cross edge $xy$ is replaced in $M_1$ with an arbitrary edge inside $V_i$, the set $M_1$ remains an induced matching. As $M_2$ is also an induced matching, it covers at most two vertices from $V_i$. We therefore have a vertex $x' \in V_i$ such that $(M \setminus \{xy\}) \cup \{xx'\}$ is an II-matching and $|M|=|M'|$. Repeating this procedure for all cross edges in $M_1 \cup M_2$, we obtain an II-matching without cross edges. Again, we may refer to the property that if $V_i$ contains an edge from $M_p$, for $p \in [2]$, then $V(M_p) \cap V_{i'}= \emptyset$ whenever $v_iv_{i'} \in E(F)$. We may conclude that $|M_p| \leq \alpha(F)$ and hence, $|M| \leq 2\alpha(F)$. This finishes the proof for  $\mbox{ii}(G_F)= 2\alpha(F)$.

\paragraph{Conclusion.} As $G_F$ is obtained by a polynomial-time construction from $F$, and the NP-complete problem of deciding whether $\alpha(F)\geq k$ holds can be reduced to the problem of deciding about $\mbox{ii}(G_F) \geq 2k$, the latter problem is also NP-hard. It is also clear that the decision problem of  $\mbox{ii}(G) \geq \ell$ belongs to NP. Thus, we may infer that the problem is NP-complete over the class of all graphs. 
\qed

\subsection{$(4t+2)$-subdivisions} \label{sec:4t+2}
 
\begin{theorem}
 	\label{thm:4t+2}
 	For every graph $G$ and integer $k \equiv 2 \bmod 4$, 
 	\begin{displaymath}
 	\gsp(S_k(G))= 
 	\frac{k}{2}\, m(G)+ n(G)- \mbox{ii}(G).
 	\end{displaymath}
 \end{theorem}
 \proof
 Let $V(G)=\{v_1,\ldots,v_n\}$, and $k=4t+2$. Choose a $\gsp$-set $D$ of $S_k(G)$. Let $A=\overline{D}$ and $B$ an arbitrary core of $D$. By Lemma~\ref{lem:matching}, the edges $E[A,B]$ form a matching $M$ in $S_k(G)$.  Again, we have three possibilities for a super edge $P_{v_iv_j}$.
 \begin{itemize}
 	\item[$(a)$] $M$ contains at most $\frac{k}{2}-1$ edges from $P_{v_iv_j}$. Let $E_0$ denote the set of the edges $v_iv_j \in E(G)$ with this property.
 	\item[$(b)$] $M$ contains exactly $\frac{k}{2}$ edges from $P_{v_iv_j}$. The set of the corresponding  edges $v_iv_j$ in $G$ is denoted by $E_1$.
 	\item[$(c)$] $M$ contains exactly $\frac{k}{2}+1$ edges from $P_{v_iv_j}$. In this case, both $v_i(v_iv_j)_1$ and $(v_iv_j)_kv_j$ belong to $M$, and we set $v_iv_j \in E_{2}$.	 
 \end{itemize}
 By definitions, the following inequality holds:
 \begin{eqnarray} \label{ineq1:4t+2}
 |M|\leq \frac{k}{2}\, m(G) + |E_2|-|E_0|. 
 \end{eqnarray}
 
 For every edge $v_iv_j \in E_2$ of $G$, the super edge $P_{v_iv_j}$ must contain the edges $v_i(v_iv_j)_1, (v_iv_j)_2(v_iv_j)_3, \dots ,  (v_iv_j)_{4t}v_j$  from $M$. Moreover, if $v_i \in A$ in $S_k(G)$ then, by Lemma~\ref{lem:matching}, $E[A,B]$ is a matching in $S_k(G)$ and 
 $$B \cap V(P_{v_iv_j})= \{(v_iv_j)_{1}, (v_iv_j)_{2}, (v_iv_j)_{5}, (v_iv_j)_{6}, \dots, (v_iv_j)_{4t+1}, (v_iv_j)_{4t+2}\}.$$ 
  The remaining subdivision vertices and $v_j$ then belong to $A$. Therefore, $v_i \in A$ implies $v_j \in A$ and, similarly, $v_i \in B$ implies $v_j \in B$ if $v_iv_j \in E_2$. We may therefore partition $E_2$ into 
  $$E_A=\{v_iv_j: v_iv_j \in E_2 \mbox{ and } v_i, v_j \in A\} \quad \mbox{and} \quad E_B=\{v_iv_j: v_iv_j \in E_2 \mbox{ and } v_i, v_j \in B\}.$$
 
 Suppose now that $E_A$ contains two edges $v_iv_j$ and  $v_pv_q$ and there exists an edge $v_iv_q \in E(G)$. We will say that $v_iv_q$ is an \emph{$A$-critical edge}. As $v_i$ and $v_q$ are already covered by $M$-edges, $v_iv_q \notin E_2$. By the same reason, if $v_iv_q \in E_1$, then the $k/2$ edges in $M \cap E(P_{v_iv_q})$ are $(v_iv_q)_1(v_iv_q)_2, \dots , (v_iv_q)_{4t+1}(v_iv_q)_{4t+2}. $ 
 Referring to Lemma~\ref{lem:matching} again, $v_i \in A$ implies $(v_iv_q)_{1} \in A$, $(v_iv_q)_{2}, (v_iv_q)_{3} \in B,  \dots, (v_iv_q)_{4t+2} \in B$, and also that $v_q \in B$. As $v_q \in A$ was supposed, it is a contradiction. We may infer that $v_iv_q \in E_0$ holds for every $A$-critical edge. The same is true for the set $B$; that is if $v_iv_j, v_pv_q \in E_B$ and a $B$-critical edge $v_iv_q$ is present in $G$, then $v_iv_q \in E_0$.
 
 We now prove that the maximum for $|E_2|-|E_0|$  can be attained without the presence of $A$- and $B$-critical edges. Indeed, let $v_iv_q$ be an $A$-critical edge such that $v_iv_j, v_pv_q \in E_A$. We may perform the following alteration in $M$ without decreasing $|E_2|-|E_0|$. First we remove the edges $M\cap E(P_{v_i,v_j})$ from $M$ and replace them by the complement edge set $E(P_{v_1,v_j}) \setminus M$. Then, we also remove the at most $\frac{k}{2}-1$ edges $M \cap E(P_{v_iv_q})$ and replace them by the $\frac{k}{2}$ edges $(v_iv_q)_{1}(v_iv_q)_{2}, \dots , (v_iv_q)_{4t+1}(v_iv_q)_{4t+2}$. After these changes we update the sets $A$ and $B$ along $P_{v_iv_j}$ and $P_{v_iv_q}$ such that, for $x \in \{j,q\}$ and $s \in [4t+2]$, a subdivision vertex $(v_iv_x)_s$ belongs to $B$ if $s \bmod 4 \in \{0,1\}$, otherwise it is put into $A$. The vertices $v_i$, $v_j$ are not in $A$ anymore, but we still have $v_q \in A$. It can be checked that $E[A,B]$ remains a matching and, by Lemma~\ref{lem:matching}, $D=\overline{A}$ is a super dominating of $S_k(G)$.  
 By this change, both $v_iv_j$ and $v_iv_q$ are moved to $E_1$, the number of $A$-critical edges is decreased, while $|E_2|-|E_0|$ remains the same.
 By the symmetry of the roles of the sets $A$ and $B$, if a $B$-critical edge exists, we may do the analogous changes. 
 
 Repeating these changes while there are critical edges, we obtain a matching $M$ and a super dominating set $D$ without decreasing $|E_2|-|E_0|$. The edges in $E_2$ still form a matching as every $u \in V(G)$ is covered only one edge from $M$; and the absence of $A$- and $B$-critical edges means that both $E_A$ and $E_B$ are induced matchings in $G$. It implies that the obtained $M$ is an II-matching and hence, $|E_2| \leq \mbox{ii}(G)$. From (\ref{ineq1:4t+2}), we now obtain
  $$ |M|\leq \frac{k}{2}\, m(G) + \mbox{ii}(G),
 $$
 and we may conclude
 \begin{align} \label{ineq2:4t+2}
 \gsp(S_k(G))  & =  n(S_k(G)) - |M| \geq \frac{k}{2}\, m(G) + n(G)- \mbox{ii}(G).
  \end{align} 
 \medskip
 
To complete the proof, we show that there exists a super dominating set $D$ of the required cardinality in $S_k(G)$.  Let $M^*=M^*_1 \cup M^*_2$ be a maximum II-matching in $G$. The set $D$ is constructed by the following five rules. 

 \begin{itemize}
 	\item[$(i)$] A vertex $v_i\in V(G)$ belongs to $D$ if and only if $v_i \notin V(M_1^*)$.
 	\item[$(ii)$] If $v_iv_j \in M^*_1$, then $D$ contains the following vertices from the super edge $P_{v_iv_j}$:
 	$$ (v_iv_j)_1, (v_iv_j)_2, (v_iv_j)_5, (v_iv_j)_6, \dots, (v_iv_j)_{4t+1}, (v_iv_j)_{4t+2}. 	$$
 	\item[$(iii)$] If $v_iv_j \in M^*_2$, then $D$ contains the following vertices from $P_{v_iv_j}$:
 	$$ (v_iv_j)_3, (v_iv_j)_4, (v_iv_j)_7, (v_iv_j)_8, \dots, (v_iv_j)_{4t-1}, (v_iv_j)_{4t}. 	$$
 	\item[$(iv)$] If $v_iv_j \notin M^*$ and $v_i \in V(M_1^*)$, then $v_j \notin V(M_1^*)$. In this case, we put the following subdivision vertices into $D$:
 	$$ (v_iv_j)_2, (v_iv_j)_3, (v_iv_j)_6, (v_iv_j)_7, \dots, (v_iv_j)_{4t-2}, (v_iv_j)_{4t-1}, (v_iv_j)_{4t+2}. 	$$
 	\item[$(v)$] If $v_iv_j \notin M^*$ and $v_i \in V(M_2^*)$, $v_j \notin V(M^*)$, and also if $v_i, v_j \notin V(M^*)$, we put the following subdivision vertices into $D$:
 	$$ (v_iv_j)_1, (v_iv_j)_4, (v_iv_j)_5,  \dots, (v_iv_j)_{4t}, (v_iv_j)_{4t+1}. 	$$
  \end{itemize}

In step $(i)$ we put $n(G)-2|M_1^*|$ non-subdivision vertices into $D$; in step $(ii)$ we consider $|M_1^*|$ super edges and put $\frac{k+2}{2}$ subdivision vertices from each into $D$; in $(iii)$ $|M_2^*|$ super edges are considered and put $|M_2^*|\, \frac{k+2}{2}$ vertices into $D$. For the remaining edges $v_iv_j$ of $G$ step $(iv)$ or $(v)$ is applied. In either case, $D$ contains exactly $\frac{k}{2}$ subdvision vertices from $P_{v_iv_j}$. The size of $D$ is therefore
 \begin{align*}
 |D|  & =  (n(G)-2|M_1^*|)+ |M_1^*|\left(\frac{k}{2}+1\right) +|M_2^*| \left(\frac{k}{2}-1\right) + (m(G)- |M_1^*|-|M_2^*|)\, \frac{k}{2}\\
  & = \frac{k}{2}\, m(G)   + n(G) - (|M_1^*|+|M_2^*|)\\
  & = \frac{k}{2}\, m(G)  + n(G) - \mbox{ii}(G).
 \end{align*} 
 It is straightforward to check that $D$ is a super dominating set in $S_k(G)$. We notice that if a vertex $v_i \in V(G)$ does not belong to $D$, then there exists an edge $v_iv_j \in M_1^*$ and, by $(ii)$, the subdivision vertex $(v_iv_j)_1$ super dominates $v_i$. If $v_i \in V(M_2^*)$ such that $v_iv_j \in M_2^*$, then $v_i \in D$ and $(v_iv_j)_1$ is $v_i$'s  only neighbor which is not in $D$. Then, $v_i$ super dominates $(v_iv_j)_1$ (that is the same as $(v_jv_i)_{4t+2}$.) For a super edge $P_{v_iv_j}$ considered in $(v)$, both ends $v_i$ and $v_j$ belong to $D$ and hence, $(v_iv_j)_1$ super dominates $(v_iv_j)_2$. 
 
 Since the constructed set $D$ is a super dominating set, we may conclude  
 $$\gsp(S_k(G)) \leq |D|  =\frac{k}{2}\, m(G) + n(G) -\mbox{ii}(G)$$
  which, together with (\ref{ineq2:4t+2}), complete  the proof of the theorem. \qed 
  \medskip
  
  As a consequence of Proposition~\ref{prop:ii} and Theorem~\ref{thm:4t+2}, we obtain the following result.
  
  \begin{theorem} \label{comp:4t+2}
  	Over the class of $k$-subdivision graphs, it is NP-complete to decide whether $\gsp(F) \leq \ell$ holds, if $\ell$ is part of the input and  $k$ is a fixed integer with  $k \equiv 2 \bmod 4$. 
  \end{theorem}

\section{Conclusions}
We conclude the paper by summarizing our main results on the computational complexity of the following problem.

\bigskip
\medskip

\noindent\fbox{%
	\parbox{\textwidth}{%
		\noindent \textbf{SUPER DOMINATION (S-DOM) PROBLEM}
		\begin{itemize}
			\item[\empty]\textit{Instance:} A simple undirected graph $G$ and an integer $\ell$.
			\item[\empty]\textit{Question:} Does $\gsp(G) \leq \ell$ hold?
		\end{itemize}
	}%
}
\medskip

\noindent
By Theorems~\ref{thm:tree}, \ref{thm:bipartite}, \ref{comp:4t+3}, \ref{comp:4t+1}, \ref{comp:4t}, and \ref{comp:4t+2}, we can conclude the following summary.
\medskip
\medskip

\noindent\fbox{%
	\parbox{\textwidth}{%
\begin{itemize}
\item The S-DOM problem is NP-complete over the following graph classes:
\begin{itemize}
\item[(A)] Bipartite graphs of girth at least $8$;
\item[(B)] Class of $k$-subdivision graphs if $k$ is a fixed even integer.
\end{itemize}
\item The S-DOM problem can be solved in polynomial time over the following graph classes:
\begin{itemize}
	\item[(C)] Trees;
	\item[(D)] Class of $k$-subdivision graphs if $k$ is an odd integer.
\end{itemize}
\end{itemize}
}%
}
\medskip

\noindent
Classes (A), (C), and (D) are subclasses of bipartite graphs. However, to get a better picture of the computational complexity of the S-DOM problem on the class of bipartite graphs, we propose the following problems.
\begin{problem}
Find a subclass ${\cal C}$ of (A) so that the S-DOM problem remains NP-complete on ${\cal C}$.
\end{problem}
\begin{problem}
	Find further subclasses of bipartite graphs over which the S-DOM problem can be solved in polynomial time. 
\end{problem}
The ultimate goal we set is a complete characterization:
\begin{problem}
	Characterize the subclasses of bipartite graphs over which the S-DOM problem remains NP-complete.
\end{problem}

\section*{Acknowledgments} 

Nima Ghanbari would like to thank the Research Council of Norway and Department of Informatics, University of Bergen for their support. Csilla Bujt\'as and Sandi Klav\v zar were supported by the Slovenian Research Agency (ARRS) under the grants P1-0297,  J1-2452, and N1-0285.

\section*{Declaration of interests}
 
The authors declare that they have no conflict of interest. 
\section*{Data availability}
 
Our manuscript has no associated data.

\end{document}